\documentclass{article} \usepackage{amssymb,amsfonts}
\title{\huge\bf On zeros of polynomials and allied functions satisfying second order differential equations}
      
\author{   {\bf Ilia Krasikov}\\
		 Brunel University\\
			    Department of Mathematical Sciences\\
				       Uxbridge UB8 3PH United Kingdom\\
						  e-mail: mastiik@brunel.ac.uk
							 }
							  
							  \date{}

\newcommand{\LP}{\mbox{$\cal L - \cal P $}}

\newcommand{\QED}{\hfill$\Box$}

\newcommand{\be}{\begin{equation}}
\newcommand{\ee}{\end{equation}}

\newtheorem{lemma}{Lemma}
\newtheorem{theorem}{Theorem}

\mathchardef\inn="3232
\renewcommand{\in}{\mbox{$\,\inn\,$}}

\begin{document}
%\title[Polynomials ]{On zeros of polynomials}
    
%\author[I. Krasikov]{Ilia Krasikov}
     
%\address{   Department of Mathematical Sciences,
%            Brunel University,
%            Uxbridge UB8 3PH United Kingdom}
%\email{mastiik@brunel.ac.uk}
			      
%\subjclass{33C45}
			       
%\date{}
%%%%%%%%%%%%%%%%%%%%%%%%%%%%%%%%%%%%%%%%%%%%%
				
%-------------------- Title Page ----------------------------------

%------------------------------------------------------------------------%
%                                                                        %
%         1. INTRODUCTION                                                %
%                                                                        %
%------------------------------------------------------------------------%
\maketitle
\vspace*{4ex}
 
\begin{center}
{\bf Abstract}
\end{center}
  
\begin{quote}
\baselineskip3ex
We shall give bounds on the spacing of zeros
of certain functions 
belonging to the  Laguerre-P\'{o}lya class and 
satisfying a second order
linear differential equation.
As a corollary we establish new sharp inequalities on the extreme zeros of
the Hermite, Laguerre and Jacobi polynomials, which are uniform in all the parameters.
\end{quote}
\vspace{16ex}
 
\noindent

%================= Footnotes =================
%%%%%%%%%%%%%%%%%%%%%%%%%%%%%%%%%%%%%%%%%%%%%%%%%%%%%%%%%%%%%%%%%%%%%%%%%%%%%%%%%%%%%%%%%%%%%%%%%%%%%%5
\section{Introduction}
The aim of this paper is to establish new sharp inequalities on the extreme zeros of
the classical orthogonal polynomials which are uniform in all the parameters.
We shall use a modification of the method suggested in \cite{k1}. 
In fact, our result is more general, and deals with the solutions
of the second order differential equation with variable coefficients
\begin{equation}
\label{difeq}
f''-2a f'+b f=0,
\end{equation}
belonging to the  Laguerre-P\'{o}lya class  \LP.
The Laguerre-P\'{o}lya class consists of real polynomials having only real zeros and real entire functions
having a representation of the form
\begin{equation}
\label{lagpolc}
c x^m e^{-\alpha x^2+\beta x} \prod_{i=1}^{\infty} (1-\frac{x}{x_i}) e^{x/{x_i}}  \, ,
\end{equation}
where $c, \beta, x_i$ are real, $\alpha \ge 0$, $m$ is a nonnegative integer and
$\sum x_i^{-2} < \infty.$
The well-known inequality of Laguerre (\cite{lagu}, p.171) states
that 
\begin{equation}
\label{laguq}
U(f)=f'^2-f f'' = - \, \left( \frac{f'}{f} \right)'  f^2 \ge 0,
\end{equation}
for any $f \in$\LP. 
\\
We shall show that this simple inequality readily yields quite accurate bounds on the spacing of
zeros of $f.$ Our main result, Theorems \ref{le1} and  \ref{nov} below, will be proved  
in the next section. Then we use Theorem \ref{le1} to derive bounds on the extreme zeros of the
Bessel function, generalized Hermite, Laguerre and Jacobi polynomials. 
As Theorem \ref{le1} does not provide any information
concerning  the precision of the inequalities,
we will compare our results with the existing bounds, 
whenever the corresponding asymptotic or
inequalities are known.
It seems that the typical situation, at least for 
the classical orthogonal polynomials with parameters growing not faster than the degree, is as following.
The extreme zero, say the largest one  $x_M ,$ of $f(x),$ has the asymptotic expansion
$x_M \thicksim P-3^{-1/3} i_{11} Q k^{-2/3}+...,$ where $k=deg (f)$
in the polynomial case. Here
$3^{-1/3} \,  i_{11}=1.855757...,$ and $i_{11}$ denotes
the least positive zero of the Airy function
$$A(t)=\frac{\pi}{3} \, \sqrt{\frac{t}{3}} \, \left(
J_{-1/3} (s)+J_{1/3} (s)
\right) , $$
where $s=2 (\frac{t}{3})^{3/2}, $ \cite{ab}.
On the other hand, in a few cases we were able to compare with, (\ref{eqglM}) gives
$x_M < P- \frac{3}{2} Q k^{-2/3},$ with the same values of $P$ and $Q.$
Unfortunately, only the first term
of the asymptotics is known when $f$ contains parameters
which can vary with $k.$ (see e.g. \cite{det, fal, ismail, assch} and the references therein). It is tempting to conjecture
that even in this case our second term is sill close to the correct value.

\section{Main Theorem}
Given an $f\in \LP$ satisfying (\ref{difeq})
we shall introduce two functions, the logarithmic derivative $t(x)=\frac{f'(x)}{f(x)}$
and the discriminant $\Delta (x)=b(x)-a^2(x).$
If $f(x)$ is a polynomial of degree $m \ge k$ 
with $k$ distinct zeros $x_1 <x_2 ...<x_k ,$ counting without multiplicity,
put formally $x_0=-\infty ,x_{k+1}=\infty .$
Then, by (\ref{laguq}), $t(x)$ consists of $k+1$
decreasing branches $B_0,B_1,..., B_k ,$
where $B_i$ is defined on $(x_{i-1}, x_i).$ 
The same notation with an obvious modification will be used for the case of entire functions,
when the sequences $\{ x_i \}$ and $\{ B_i \}$
are one or both side infinite.
\begin{theorem}
\label{le1}
Let  $f \in \LP$ satisfy (\ref{difeq}) and
suppose that $a(x)$ intersects a branch $B_i$ for $x=c_i .$
Let $J$ be the region defined by $\Delta (x) > 0.$
Then $c_i \in J,$ and moreover,
\begin{equation}
\label{eq11}
 x_i< c_i-\frac{1}{\sqrt{\Delta (c_i)}} \, , 
\end{equation}
\begin{equation}
\label{eq12}
x_{i+1}> c_i+\frac{1}{\sqrt{\Delta (c_i)}} \, . 
\end{equation}
In particular, if $a(x)$ intersects either $B_0$ or $B_k$ then 
\begin{equation}
\label{eqglm}
x_1 > \min_{x \in J} \{ x+\frac{1}{\sqrt{\Delta (x)}}\}, \; \; \;
\end{equation}
\begin{equation}
\label{eqglM}
x_k < \max_{x \in J} \{ x-\frac{1}{\sqrt{\Delta (x)}}\} \, ,
\end{equation}
respectively.
\end{theorem}
\begin{proof}
Let $x_j$ be a zero of $f$, consider $g(x)=f(x)/(x-x_j).$
Using $f''=2a f'-b f,$ we get
$$
0 \le (x-x_j)^4 U(g)=(x-x_j)^4 ({g'}^2-g g'')=(x-x_j)^2 {f'}^2-2(x-x_j)^2 a f f'+((x-x_j)^2 b-1)f^2.
$$
Dividing by $f^2$
yields
$$(x-x_j)^2 ( t^2(x)-a(x) t(x)+b(x)) \ge 1.$$
Therefore, for any $c$ being a solution of $t(x)=a(x),$
one has
\begin{equation}
\label{vaj}
(c-x_j)^2(b(c )-a^2(c))\ge 1.
\end{equation}
This implies $b(c)>a^2(c).$
Now, if $c=c_i,$ that is $x_i <c_i <x_{i+1},$ choosing $j=i$ or $i+1$ in (\ref{vaj}),
we conclude 
$x_i < c_i-\frac{1}{\sqrt{\Delta (c_i)}}$,
$x_{i+1} > c_i+\frac{1}{\sqrt{\Delta (c_i)}}, $
thus proving (\ref{eq11}),(\ref{eq12}).
Finally (\ref{eqglm}),(\ref{eqglM}) follow by
$$x_1 > c_0+\frac{1}{\sqrt{\Delta (c_0)}} \ge \min_{x \in J} \{ x+\frac{1}{\sqrt{\Delta (x)}}\}\, ,$$
$$x_k < c_k-\frac{1}{\sqrt{\Delta (c_k)}} \le \max_{x \in J} \{ x-\frac{1}{\sqrt{\Delta (x)}}\} \, .$$
\QED
\end{proof}

To understand what type of bounds maybe derived from (\ref{eqglm}) and (\ref{eqglM}),
suppose that $\Delta (x)$ has precisely two zeros ${\bf y_1}<{\bf y_2}.$ Let the minimum
in (\ref{eqglm}) be attained for $x = {\bf y_1}+\epsilon .$
On omitting the higher order terms we have
$\Delta ({\bf y_1}+\epsilon ) \approx \epsilon \, \Delta' ({\bf y_1}) ,$ and so
$$ x_1 > \min  \{ x+\frac{1}{\sqrt{\Delta (x)}}\} \approx {\bf y_1}+ 
\min_{\epsilon >0} \{ \epsilon +\frac{1}{ \sqrt{ \epsilon \, \Delta' ({\bf y_1}) }} \} = {\bf y_1}+3  (4 \Delta' ({\bf y_1}))^{-1/3} \,.
$$
Similarly, we get
$x_k \lessapprox {\bf y_2}-3  (4 \Delta' ({\bf y_2}))^{-1/3} \,.$

Therefore one could expect, say, for the least zero, that there are constants $A, B, C$ such that 
$A < \frac{x_2-x_1}{x_1-{\bf y_1}} <B$ and 
$x_1-{\bf y_1} =C (\Delta' ({\bf y_1}))^{-1/3}.$
\\
Theorem \ref{le1} also implies that
\begin{equation}
\label{eq22}
x_{i+1}-x_i >  \min_{x_i <x <x_{i+1}} \frac{2}{\sqrt{\Delta (x)}} \, .
\end{equation}
Slightly stronger result can be proved if we consider
$U(g)$ with $g(x)=\frac{f(x)}{(x-x_{i})(x-x_j )}.$

\begin{theorem}
\label{rast}
Let  $f \in \LP$ satisfy (\ref{difeq}), have only simple zeros, and
suppose that $a(x)$ intersects $t(x)$ between zeros $x_i < x_j$ of $f.$ Then
\begin{equation}
\label{urraz1}
(x_{i}-x_j)^2 \ge  \min_{x_i <x <x_{j}} \frac{8}{\Delta (x)} .
\end{equation}
\end{theorem}
\begin{proof}
Set $g(x)=\frac{f(x)}{(x-x_{i})(x-x_j )},$ then using (\ref{difeq})
to eliminate higher derivatives of $f$, we obtain
$$U(g)=
\frac{(x-x_{i})^2 (x-x_j )^2 \left( f'^2(x)-2a(x)f'(x)f(x)+b(x)f^2(x) \right) -
\left( (x-x_{i})^2+(x-x_j )^2 \right) f^2(x)}{(x-x_{i})^4 (x-x_j )^4}.
$$
Consider this expression at a point $c,$
$a(c)=t(c),$ and $x_i <c <x_j.$
Using $(c-x_{i})^2+(c-x_j )^2 \ge 2(c-x_{i})(x_j-c),$ 
we obtain
$$
\Delta (c) \ge \frac{2}{(c-x_{i})(x_j-c)} \ge \frac{8}{(x_j-x_i)^2},
$$
yielding (\ref{urraz1}).
\QED
\end{proof}

We expect that in many cases (\ref{urraz1}) is of the correct order besides the factor $8.$
For example, for $f= \sin{ x}$ ( $f''+f=0$) it gives $x_{i+1}-x_i \ge 2 \sqrt{2},$ instead of $\pi .$
It is easy also to check that for the Chebyshev polynomials $T_k$
the true answer can be at most $\sqrt{\frac{\pi^2-1}{2}} \approx 2.1$
times greater than that given by (\ref{urraz1}).

Theorem \ref{le1} has one shortcoming. The graph of $t(x)$ consists of
cotangent-shape branches in the middle and hyperbolic branches at the ends.
Whenever the condition on intersection of $t(x)$ with $a(x)$ at a cotangent-shape branch
is almost automatically fulfilled, the intersection with the uttermost hyperbolic branches
is not obvious (fortunately, for classical orthogonal 
polynomials $a(x)$ does intersects all the branches).
We can get rid of the intersection
conditions
if we restrict the class of functions $f$ and 
assume that $a(x)$ and $b(x)$ are
sufficiently smooth in a vicinity of zero $x_i$ of $f.$
Namely, we consider entire functions of order less than 2 with only distinct real zeros.
By Hadamard's factorization theorem (see e.g. \cite{tich}),
such functions are either polynomials or have a canonical product representation 
\begin{equation}
f(x)=c x^m  e^{\beta x} \, \prod_{i=1}^{\infty} (1-\frac{x}{x_i}) e^{x/{x_i}} \, ,
\end{equation}
where $c, \beta, x_i$ are real, $m=0$ or $1,$ and
$\sum x_i^{-2} < \infty,$
We denote this class by $\LP (I).$
It is well known (the result usually attributed to Laguerre, see e.g. \cite{tich}, p. 266,
the polynomial case is given in \cite{ps}, chapter 5, problems 62, 63),
that $f \in \LP (I)$ implies 
$$f- \lambda f' = - \lambda e^{x/\lambda} \left( e^{-x/\lambda} f(x) \right)' \in \LP (I).$$ 
Iterating this yields
\begin{lemma}
\label{polya}
Let  $f \in \LP (I),$ then for any polynomial
$$p(x)= \prod_{i=1}^n (x- \lambda_i )=\sum_{i=0}^n q_i x^i ,$$
$\lambda_1 , ..., \lambda_n ,$ are real,
the function $g=\sum_{i=0}^n q_i f^{(i)} \in \LP (I),$ 
and thus $U(g) \ge 0.$
\end{lemma}
\begin{theorem}
\label{nov}
Let  $f \in \LP (I)$ satisfy (\ref{difeq}).
\\
(i) If $a(x)$ and $b(x)$ are differentiable in a vicinity of a zero $x_i $ of $f,$
then
\begin{equation}
\label{eqmin0}
\Delta(x_i)-2a'(x_i) \ge 0.
\end{equation}
(ii) If $a(x)$ and $b(x)$ are two times  differentiable in a vicinity of a zero $x_i $
of $f,$
then for $x=x_i ,$
$$\min_{\lambda} \{ 
( b^2  -8 a^2 a'-4b a'+4a'^2 +4a b'-4a a'' ) \lambda^4 -
4 ( a b -4a a'+b'-a'' ) \lambda^3+
2 ( 2a^2+b-2a' ) \lambda^2 -
$$
\begin{equation}
\label{eqmin}
4 a \lambda+1
\} \ge 0.
\end{equation}
\end{theorem}
\begin{proof}
Consider $U(g),$ where $g=f- \lambda f'.$
We have for $x=x_i ,$
$$ U(g)= ((b (x_i) -2 a'(x_i)) \lambda^2 -2 \lambda a(x_i) +1)f'^2(x_i) \ge 0,
$$
for any real $\lambda .$
Obviously, $b (x_i) -2 a'(x_i)$ must be positive, hence we can choose $\lambda =\frac{a(x_i)}{b (x_i) -2 a'(x_i)}.$
This yields (\ref{eqmin0}).
\\
To prove (\ref{eqmin}) we apply the  previous lemma with $n=2$ and $\lambda_1=\lambda_2=\lambda .$
Then $g(x)=f(x)- 2 \lambda f'(x)+ \lambda^2 f''(x),$ and the result follows
by calculating $U(g) \ge 0$ for $x=x_i.$
\QED
\end{proof}

It seems that (\ref{eqmin})  leads to the same type of bounds for the extreme zeros
as (\ref{eqglm}) and (\ref{eqglM}). 
Moreover, numerical evidences suggest that one can reach the true value
of the first two terms in the asymptotic expansion of the extreme zeros as a limiting case.
For we choose $p(x)=(x-\lambda)^n$  in Lemma \ref{polya},
e.g. $g= \sum_{i=0}^n (-\lambda)^i {n \choose i} f^{(i)},$
and consider the inequality $U(g)  >0$ at a zero $x_i .$ 
We will illustrate this for the case of Hermite polynomial.

\section{Applications}
In this section we shall use (\ref{eqglm}) and (\ref{eqglM})
to give new bounds on the extreme zeros of classical orthogonal polynomials.
We refer to \cite{det, szego}, and the references therein for the known asymptotic results,
and to \cite{ab, szego} for all formulae concerning special functions which are used in the sequel.

To gain some impression about the sharpness of the inequalities of Theorem \ref{le1}
we start with the Bessel function. 
In this case extremely precise bounds, far better than can be obtained by our method,  are known
\cite{heth, lw, lu}.
\\
{\bf Bessel functions} $J_\nu (x)$
can be defined by
the following product representation 
$$J_{\nu}(x)= \frac{x^{\nu}}{2^{\nu} \; \Gamma (\nu+1)} \, \prod_{i=1}^\infty \left( 1- \frac{x^2}{j^2_{\nu ,i}} \right) ,$$
where $j_{\nu ,1}< j_{\nu ,2}<...,$ are the positive zeros of $J_{\nu}(x).$
Thus, $u= u(x)=x^{- \nu} J_{\nu}(x)$ is an entire function and moreover
$u \in $\LP .
It can be shown directly, using
$$x^2 \, J''_{\nu}(x)+x \, J'_{\nu}(x)+(x^2- \nu^2 )J_{\nu}(x)=0, $$
that
$$x u'' +(2 \nu +1)u'+x u=0.$$
The corresponding calculations are very simple. 
\begin{theorem}
\label{besth}
\begin{equation}
\label{besoz}
j_{\nu ,1} > \frac{((2\nu+1)^{2/3}+2^{2/3})^{3/2}}{2}
\end{equation}
provided $\nu >- \frac{1}{2}.$
\end{theorem}
\begin{proof}
Assuming $\nu >- \frac{1}{2}$
and $x >0,$
one readily sees that $a(x)=- \, \frac{2 \nu +1}{2x}$ intersects all the branches $B_i$ of $t(x) $ 
for $i \ge 1.$ Using the power series representation
$$
u(x)= \sum_{i=0}^{\infty} (-1)^i \; \frac{x^{2i}}{2^{\nu +2i}i!  \Gamma (i+ \nu +1)},
$$
we find $t(0) = 0.$ Since $t(x)$ is a decreasing function tending to $-\infty$ for $x \rightarrow j_{1,\nu}^{(-)} \, ,$
it follows that $a(x)$ intersects the branch $B_0$ as well.
Now, the condition $\Delta (x) >0$ yields
$j_{\nu ,1} >\nu +\frac{1}{2},$ and moreover
$$j_{\nu ,1} > \min_{0<x<\nu +\frac{1}{2} } \left\{ x+\frac{1}{\sqrt{\Delta (x)}} \, \right\}  =
\min_{x >0} \{ x+\frac{2x}{\sqrt{4x^2-(2 \nu +1)^2}}\}= \frac{((2\nu+1)^{2/3}+2^{2/3})^{3/2}}{2} \,,$$
where the minimum is attained for
$$x= \frac{(2 \nu +1)^{2/3} \; \sqrt{(2 \nu +1)^{2/3}+ 2^{2/3}}}{2} \,.$$
\QED
\end{proof}

The bound given by (\ref{besoz}) is $ \nu+\frac{3}{2} \nu^{1/3}+O(1) $
for $k \rightarrow \infty .$
On the other hand, it is known \cite{heth} that the first two  (in fact three, see \cite{lw, lu})
terms of the asymptotic expansion of $j_{\nu ,1}$ provide a lower bound for it.
Namely, for $\nu >0,$  $j_{\nu ,1} >  \nu+3^{-1/3} \,  i_{11} \,  \nu^{1/3} .$
Thus, (\ref{besoz}) gives the correct answer up to the value of the constant at the second term,
$1.5$ instead of $1.855757... \,$. 

{\bf Generalized Hermite polynomials} $H_k^{\mu} (x)$ are polynomials orthogonal on $(-\infty ,\infty )$ for $\mu > - \frac{1}{2},$
with respect to the weight function $|x|^{2 \mu} e^{-x^2}.$ The corresponding ODE is
$$
u''-2(x-\mu x^{-1}) u'+(2k - \theta_k x^{-2} )u=0, \; \; u=H_k^{\mu} (x),
$$
where $\theta_{2i}=0, \; \theta_{2i+1}=2 \mu .$
\\
The following result is an improvement on asymptotics given in \cite{det, fal}.
\begin{theorem}
\label{genher}
Let $x_m$ and $x_M$ be the least and the largest positive zero of $H_k^{\mu} (x)$ respectively, $\mu > - \frac{1}{2}.$
Then
\begin{equation}
\label{ozhrm}
x_m > \sqrt{k+\mu - r} +\frac{3}{2} \left( 
\frac{k+\mu -r }{4 r^2} \right)^{1/6} \, ,
\end{equation}
\begin{equation}
\label{ozhrM}
x_M< 
\sqrt{k+\mu + r} -\frac{3}{2} \left(
\frac{k+\mu + r }{4 r^2} \right)^{1/6} \, ,
\end{equation}
where $r=\sqrt{k^2+2k \mu-\theta_k } \, .$
\end{theorem}
\begin{proof}
The zeros of $H_k^{\mu} (x)$ are symmetric with respect to the origin, hence we may assume $x>0.$
Since $a(x)=x-\mu x^{-1},$ is a continuous increasing function 
for $\mu \ge 0,$ and positive tending to $\infty$ for $ \mu <0,$ 
in this region,
it intersects all the branches $B_i$
corresponding to the positive zeros of $H_k^{\mu} (x).$
The discriminant
$$ \Delta(x)=\left( 2k x^2- \theta_k -(x^2- \mu )^2 \right) x^{-2},$$
has two positive roots ${\bf y_{1,2}}=\sqrt{k+\mu \pm r}.$
Solving $\Delta (x) >0,$
by Theorem \ref{le1}, we obtain,
${\bf y_1} <x_m <x_M<{\bf y_2}.$
Moreover,
\begin{equation}
\label{ozhrm1}
x_m > \min_{{\bf y_1} <x<{\bf y_2}} \{ x+ \frac{1}{\sqrt{\Delta (x)}} \},
\end{equation}
\begin{equation}
\label{ozhrM1}
x_M < \max_{{\bf y_1} <x<{\bf y_2}} \{ x- \frac{1}{\sqrt{\Delta (x)}} \}.
\end{equation}
We shall prove here (\ref{ozhrm}), the proof of (\ref{ozhrM}) is similar.
Suppose that the minimum in (\ref{ozhrm1}) is attained for $x={\bf y_1}+\epsilon .$ 
Since $\Delta ({\bf y_1})=0,$ and $\Delta''(x)= \frac{-2(x^4+3 \mu^2 +3 \theta )}{x^4} <0,$
we have
$\Delta ({\bf y_1}+\epsilon )<   \epsilon \Delta' ({\bf y_1}) ,$
and so
$$x_m > \min_{{\bf y_1} <x<{\bf y_2}} \{ x+ \frac{1}{\sqrt{\Delta (x)}} \} > 
{\bf y_1}+\epsilon +\frac{1}{\sqrt{  \epsilon  \, \Delta' ({\bf y_1} )}} \ge
{\bf y_1}+\min_{\epsilon >0} \{ \epsilon + \frac{1}{\sqrt{ \epsilon \, \Delta' ({\bf y_1} )}} \} = $$
$$
{\bf y_1}+ \frac{3}{2^{2/3}} \left( \Delta' ({\bf y_1}) \right)^{-1/3}.
$$
Calculations yield
$\Delta' ({\bf y_1})= 4r  \,(k+\mu-r )^{-1/2} \, , $
and the result follows.
\QED
\end{proof}

Now we suppose that $k$ is large and consider the asymptotics corresponding to (\ref{ozhrm}) and (\ref{ozhrM}).
If $\mu $ is fixed then
$$ x_m > \frac{2 \sqrt{\mu^2 +\theta } +3 ( \mu^2 +\theta )^{1/6}}{2 \sqrt{2k}} \, \left( 1-O(k^{-1}) \right) .$$
$$ x_M < \sqrt{2k}-\frac{3}{2} \, (2k)^{-1/6} \left( 1- O(k^{-1/3}) \right) .$$
If $\frac{\mu}{k}=\delta $ is fixed, then
$$
x_m > (\sqrt{2 \delta +1} -1)\, \sqrt{\frac{k}{2}} + \frac{3}{2 \sqrt{2}} \, 
(\sqrt{2  \delta +1} -1)^{1/3} (2 \delta+1 )^{-1/6} k^{-1/6} \left( 1- O(k^{-1/3}) \right) \,.
$$
$$
x_M <
(\sqrt{2 \delta +1} +1)\, \sqrt{\frac{k}{2}} - \frac{3}{2 \sqrt{2}} \, 
(\sqrt{2  \delta +1} +1)^{1/3} (2 \delta+1 )^{-1/6} k^{-1/6} \left( 1- O(k^{-1/3}) \right) \,.
$$
For $\mu=0,$ a sharper result is known (see e.g. \cite{szego}, sec.6.32),
which indicates the same loss of the precision in (\ref{ozhrM}) as for the Bessel function,
\begin{equation}
\label{segh}
x_M < \sqrt{2k+1} -6^{-1/3} \, (2k+1)^{-1/6} i_{11} =\sqrt{2k+1}-1.85575 \, (2k+1)^{-1/6}.
\end{equation}
%%%%%%%%%%%%%%%%%%%%%%%%%%%%%%%%%%%%%%%%%%%%%%%%%%%%%%%%%%%%%%%%%%%%%%%%%%%%%%%%%%%%%%%%%%%%%%%%%%%%%%%%%%%%%%%%%

{\bf Laguerre Polynomials} $L_k^{(\alpha)}(x)$ are polynomials orthogonal on $[0, \infty )$ for $\alpha >-1,$
with respect to the weight function $x^{\alpha } e^{-x}.$
The corresponding ODE is
$$
u''-(1-( \alpha+1) x^{-1})u'+k x^{-1} u=0, \; \;\; u=L_k^{(\alpha)}(x).
$$
\begin{theorem}
\label{lagr}
Let $x_m$ and $x_M$ be the least and the largest zero of $L_k^{(\alpha)} (x)$ respectively, 
$\alpha >-1.$
Then
\begin{equation}
\label{lm}
x_m > r^2+3 r^{4/3} (s^2-r^2)^{-1/3}\, ,
\end{equation}
\begin{equation}
\label{lM}
x_M < s^2-3 s^{4/3} (s^2-r^2)^{-1/3}+2 \, ,
\end{equation}
where $r=\sqrt{k+\alpha+1}-\sqrt{k},$ $s=\sqrt{k+\alpha+1}+\sqrt{k} \, .$
\end{theorem}
\begin{proof}
Using the variables $r$ and $s$ we get
$a(x)=\frac{1-r s}{2x},$
$b(x)=\frac{(s-r)^2}{4x},$
$\Delta (x)=\frac{(x-r^2)(s^2-x)}{4x^2}.$
Obviously, $a(x)$ is a continuous increasing function for $\alpha >-1,$ and $x >0.$
Thus, it intersects all the branches of $t(x).$
By Theorem \ref{le1}
we have $r^2 < x_m < x_M < s^2,$
and also
\begin{equation}
\label{ozlgm}
x_m > \min_{{ r^2} <x<{s^2}} \{ x+ \frac{1}{\sqrt{\Delta (x)}} \},
\end{equation}
\begin{equation}
\label{ozlgM}
x_M < \max_{{r^2} <x<{s^2}} \{ x- \frac{1}{\sqrt{\Delta (x)}} \}.
\end{equation}
To prove (\ref{lm}), 
we assume that the minimum is attained at $x=r^2+ \epsilon$. 
We obtain
$$\Delta (r^2+ \epsilon )= \frac{\epsilon (s^2-r^2-\epsilon )}{4(r^2+\epsilon )^2} <\frac{\epsilon (s^2-r^2)}{4r^4}.$$
Therefore
$$x_m >\min_{{ r^2} <x<{s^2}} \{ x+ \frac{1}{\sqrt{\Delta (x)}} \}>r^2+\epsilon+\frac{2r^2}{\sqrt{\epsilon (s^2-r^2)}}
\ge r^2+ \min_{\epsilon >0} \{ \epsilon+\frac{2r^2}{\sqrt{\epsilon (s^2-r^2)}} \} = $$
$$
r^2+3 r^{4/3} (s^2-r^2)^{-1/3}.
$$
To prove (\ref{lM}) we set $s^2-\epsilon$ for the extremal value of $x.$
Then
$$\Delta (s^2- \epsilon )=\frac{\epsilon (s^2-r^2-\epsilon )}{4(s^2- \epsilon )} < \frac{\epsilon (s^2-r^2 )}{4(s^2- \epsilon )^2},$$
and by (\ref{ozlgM}), 
$$
x_M < s^2- \epsilon - \frac{2(s^2- \epsilon )}{\sqrt{\epsilon (s^2-r^2)}} \le
s^2- \min_{0 < \epsilon <s^2-r^2}  \{  
\epsilon +\frac{2 s^2}{\sqrt{\epsilon (s^2-r^2)}}
\} + 2 \, \max_{0 < \epsilon <s^2-r^2}   \sqrt{\frac{\epsilon}{s^2-r^2}} =$$
$$s^2-3 s^{4/3} (s^2-r^2)^{-1/3}+2 \, .$$
This completes the proof.
\QED
\end{proof}

If $\alpha $ is fixed (\ref{lm}) and (\ref{lM}) give
$$ x_m > \frac{(1+\alpha)^2+3 (1+\alpha)^{4/3}}{4k} \left( 1-O(k^{-1}) \right) \, ,$$
$$ x_M < 4 k-3 \cdot 2^{2/3} k^{1/3} \left( 1-O(k^{-1/3}) \right) \\, .$$
If $  \frac{\alpha}{k}= \delta ,$ then
$$ x_m > (\sqrt{1+\delta }-1)^2 k +\frac{3}{2^{2/3}} \, (\sqrt{1+\delta }-1)^{4/3} (1+\delta )^{-1/6} k^{1/3}\left( 1-O(k^{-1/3}) \right) \, \, ,$$
$$ x_M < 
(\sqrt{1+\delta }+1)^2 k -\frac{3}{2^{2/3}} \, (\sqrt{1+\delta }+1)^{4/3} (1+\delta )^{-1/6} k^{1/3}\left( 1-O(k^{-1/3}) \right)
\, .$$
The classical inequality  (see \cite{szego}, sec.6.32) is
$$x_M < \left( \sqrt{4k+2 \alpha+2} -6^{-1/3} \, (4k+2 \alpha+2)^{-1/6} i_{11}
\right)^2 \, ,$$
provided $| \alpha | \ge \frac{1}{4}, \; \alpha  >-1.$
This bound is sharp only if $\alpha$ is fixed. Inequalities
uniform in $\alpha$ and $k,$ giving in fact the first terms of 
(\ref{lm}) and (\ref{lM}) has been established in \cite{ismail}.
Better bounds, practically coinciding with the main term
of our inequalities, were given by \cite{ismail}.
Inequalities with the second term only slightly weaker than in (\ref{lm}),(\ref{lM}) were recently obtained by the author
using a similar but more complicated approach \cite{k1}.
%%%%%%%%%%%%%%%%%%%%%%%%%%%%%%%%%%%%%%%%%%%%%%%%%%%%%%%%%%%%%%%%%%%%%%%%%%%%%%%%%%%%%%%%%%%%%%%%%%%%%%%%%%%%%%%%%%%%%%

{\bf Jacobi Polynomials} $P_k^{( \alpha , \beta )}(x)$ are polynomials orthogonal on $[-1,1]$
for $\alpha , \beta >-1,$
with respect to the weight function $(1-x)^{\alpha}(1+x)^{\beta}.$
The corresponding ODE is
$$
u''- \frac{( \alpha +\beta +2)x+\alpha-\beta}{1-x^2} u'+\frac{k(k+\alpha + \beta +1)}{1-x^2} u=0, \; \; \; u=P_k^{( \alpha , \beta )}(x).
$$
For the known asymptotic results giving the main terms of the following Theorem \ref{jac}
one should consult \cite{det, fal, saff}. The inequalities of the same order of precision
as asymptotics seems were known only for the ultraspherical
case \cite{for, laf}.\\
To simplify the calculations we need the following claim showing that $P_k^{( \alpha , \beta )}(x)$ has a negative zero.
\begin{lemma}
Let $x_m$ and $x_M$ be the least and the largest zero of $P_k^{( \alpha , \beta )}(x)$ respectively.
Then $x_m+x_M \le 0,$ provided $\alpha \ge \beta .$
\end{lemma}
\begin{proof}
According to the Markoff theorem (see e.g. \cite{szego}, sec. 6.21), 
$\frac{\partial x_i}{\partial \alpha} <0, \; \; \frac{\partial x_i}{\partial \beta} >0,$
for any zero $x_i$ of $P_k^{( \alpha , \beta )}(x).$
As for the ultraspherical case $\alpha = \beta ,$ we obviously have $x_m+x_M=0,$
the result follows.
\QED
\end{proof}
\begin{theorem}
\label{jac}
Let $x_m$ and $x_M$ be the least and the largest zero of $P_k^{( \alpha , \beta )}(x)$ respectively, 
$\alpha \ge \beta > -1.$
Then
\begin{equation}
\label{ozjam}
x_m > \, {\bf y_1}+3 (1-{\bf y_1}^2 )^{2/3} \, (2R)^{-1/3} \, ,
\end{equation}
\begin{equation}
\label{ozjacM}
x_M < 
{\bf y_2}-3 (1-{\bf y_2}^2 )^{2/3} \, (2R)^{-1/3} +\frac{4q (s+1)}{(r^2+2s+1)^{3/2}}
\, ,
\end{equation}
where 
$$s=\alpha +\beta +1 , \; \; 
q=\alpha - \beta , \; \; r=2k+\alpha + \beta +1, 
\; \; R=\sqrt{(r^2-q^2+2s+1)(r^2-s^2)} \, ,$$
and 
$${\bf y_1}=  - \, \frac{R+q(s+1)}{r^2+2s+1} \, , \; \; 
{\bf y_2} = \frac{R -q(s+1)}{r^2+2s+1} \, .
$$
\end{theorem}
\begin{proof}
We may assume $|x| <1.$
In this interval $a(x)$ is a continuous function, and
$$\lim_{x \rightarrow -1^{(+)}} u(x)= - \infty , \; \; \; \lim_{x \rightarrow 1^{(-)}} u(x) =\infty .$$
Thus, $a(x)$ intersect all the branches of $t(x).$ 
The corresponding discriminant is
$$\Delta (x)= - \, \frac{(r^2+2s+1)x^2+2q(s+1)x+q^2+s^2-r^2}{4(1-x^2)^2} \, ,$$
with the zeros ${\bf y_1} , {\bf y_2}.$
Thus, we obtain
${\bf y_1} <x_m <x_M <  {\bf y_2} \, . $
Observe that for ${\bf y_1} <x < {\bf y_2},$ 
$$\Delta (x) = \frac{(r^2+2s+1)(x-{\bf y_1})({\bf y_2}-x)}{4(1-x^2)^2} <
\frac{(r^2+2s+1)(x-{\bf y_1})({\bf y_2}-{\bf y_1})}{4(1-x^2)^2}= \frac{R (x-{\bf y_1})}{2(1-x^2)^2}.$$
Set $x={\bf y_1}+\epsilon $ for the extreme value in (\ref{eqglm}). 
By the previous lemma $ {\bf y_1}+\epsilon <x_m \le 0,$ and thus $2 {\bf y_1}+\epsilon <0.$
This yields
$$x_m >{\bf y_1}+  \epsilon + \frac{\sqrt{2} \, (1-({\bf y_1}+
\epsilon )^2)}{\sqrt{ \epsilon R}} 
\ge {\bf y_1}+ \min_{\epsilon >0} \left\{ \epsilon + (1-{\bf y_1}^2 )\, \sqrt{\frac{2}{ \epsilon R}} \; \; \right\} - 
\frac{\sqrt{2 \epsilon } \, (2 {\bf y_1}+\epsilon )}{\sqrt{R}}  > $$
$${\bf y_1}+3 (1-{\bf y_1}^2 )^{2/3} \, (2R)^{-1/3} 
\, .  $$
Similarly, using  $x={\bf y_2}-\epsilon > {\bf y_1},$ as the extreme value in (\ref{eqglM})
and $\Delta (x) <\frac{R ({\bf y_2}-x)}{2(1-x^2)^2},$ we obtain
$$x_M < {\bf y_2} -\epsilon - \frac{\sqrt{2} \, (1-({\bf y_2}-
\epsilon )^2)}{\sqrt{ \epsilon R}} 
\le {\bf y_2}- \min_{\epsilon >0} \left\{ \epsilon + (1-{\bf y_2}^2 )\, \sqrt{\frac{2}{ \epsilon R}} \; \; \right\}-
\frac{\sqrt{2 \epsilon } \, (2 {\bf y_2}-\epsilon )}{\sqrt{R}} \le $$
$$
{\bf y_2}-3 (1-{\bf y_2}^2 )^{2/3} \, (2R)^{-1/3} - \sqrt{\frac{2}{R}} \; \; \min_{0 <\epsilon \le {\bf y_2}-{\bf y_1}}
\{ \sqrt{\epsilon}\, (2 {\bf y_2} - \epsilon ) \} .  $$
The last minimum is attained for $\epsilon ={\bf y_2}-{\bf y_1},$
and equals to $- \, \frac{4q (s+1)}{(r^2+2s+1)^{3/2}}.$
This  completes the proof.
\QED
\end{proof}

If $\alpha , \, \beta$ are fixed and $k \rightarrow \infty ,$
$$x_m >  -1+ \frac{(1+\beta )^2 +3 (1+\beta )^{4/3}}{2k^2} \left( 1-O(k^{-1}) \right) \, ,$$
$$x_M <  1- \frac{(1+\alpha )^2 +3 (1+\alpha )^{4/3}}{2k^2} \left( 1-O(k^{-1}) \right) \, .$$
In particular, for the Chebyshev polynomial $T_k (x)$ ($\alpha =- 1/2 ,$) we get
$x_M < 1-\frac{1+3 \cdot 2^{2/3}}{8k^2} +O(k^{-3}) \approx 1-0.72 k^{-2},$
instead of the correct value $x_M= \cos{ \frac{\pi}{2k}} =1-\frac{\pi^2}{8k^2}+... \approx 1-1.23 k^{-2}.$
\\
If $\frac{\alpha}{k}=A ,\; \frac{\beta}{k}=B $ are fixed and $k \rightarrow \infty ,$
$$x_m > -S+\frac{3(1-S^2)^{2/3}}{2((1+A )(1+B )(1+A+ B ))^{1/6}} \, k^{-2/3} \left( 1-O(k^{-1/3}) \right)
\, ,$$
$$x_M < T-\frac{3(1-T^2)^{2/3}}{2((1+A )(1+B )(1+A+ B ))^{1/6}}\,  k^{-2/3} \left( 1-O(k^{-1/3}) \right) \, ,$$
where
$$S= \frac{A^2-B^2 +4 \sqrt{(1+A )(1+B )(1+A+ B )}}{(A + B +2)^2}, \; \; \; 
T= \frac{B^2-A^2 +4 \sqrt{(1+A )(1+B )(1+A+ B )}}{(A + B +2)^2} \, .$$
In the ultraspherical case $\alpha =\beta$ this yields
$$x_M <\frac{\sqrt{1+2A}}{1+A} -\frac{3 A^{4/3}}{2 (1+A)^{5/3} (1+2A)^{1/6}}\, k^{-2/3} \left( 1-O(k^{-1/3}) \right) \, .$$

To demonstrate that Theorem \ref{nov} can give inequalities of the same order of precision
as Theorem \ref{le1},
consider the case of Hermite polynomials $H_k (x)=H_k^0 (x).$
From (\ref{eqmin0}) one obtains $2k-x_i^2-2 \ge 0,$ that is $|x_i| \le \sqrt{2k-2}.$
Furthermore, by (\ref{eqmin}) we have for any zero $ x_i$ of $H_k (x),$ and any $\lambda ,$
$$
\phi=4 \lambda^2 (1-2 \lambda^2 )x_i^2-4   \lambda (2k \lambda^2-4 \lambda^2+1) x_i+(2k \lambda^2-2 \lambda^2+1)^2 \ge 0.
$$
It is easy to see that the discriminant of  this expression in $\lambda ,$
that is the resultant $Result_{\lambda} \left( \phi , \frac{\partial \phi}{\partial \lambda} \right),$
must vanish
for the extremal value of $x$.
We have for the resultant
$$
2^{14} \left( 2x^2-(k-1)^2 \right) \left( 4x^6-24(k-1)x^4+(48k^2-96k+75)x^2-32(k-1)^2 \right) .
$$
To obtain the answer in a closed form one can use the substitution
$k= \frac{m^{12}-4m^6-1}{4m^6},$ that is $m=(2k+2+\sqrt{4k^2+8k+5} \, )^{1/6}.$
This yields ( using $x < \sqrt{2k}$),
$$x \le \frac{(m^4-1)^{3/2}}{\sqrt{2} \, m^3} = \sqrt{2k}-3 \cdot 2^{-11/6} k^{-1/6}+ O(k^{-1/2}),$$
for
$$\lambda =\frac{\sqrt{2} \, m^3}{(m^4-m^2-1) \sqrt{m^4-1}}.$$
This is only slightly weaker than the bound given by (\ref{ozhrm}).
\\
If one chooses $g=\sum_{i=0}^n (- \lambda )^i {n \choose i} f^{(i)},$ for $n >2,$ the 
expression for $U(g)$ contains polynomials of high degree and becomes rather
complicated. 
But for the Hermite polynomials and a few small values of $n$ 
the asymptotic for the bounds given by $U(g) >0,$
can be obtained rather easily. We have applied the following procedure
(we used $Mathematika$ for calculations).
Given an $n$, set $k=m^6/2, x=m^3-c(n)/m, \lambda =m^{-3}+d(n) m^{-5},$
and consider the coefficient of $U(g)$ at the greatest power of $m.$
This is a polynomial, $\phi$ say, in variables $c$ and $d.$ Moreover,
the discriminant of $\phi$ in $d$ vanishes at the optimal value of $d,$
thus enabling one to exclude $d.$ It is left to find possible values of $c$ as the roots of the
obtained algebraic equation
and to select an appropriate one corresponding to the real value of $d.$
The results of these calculations indicate that $c$ strictly increases with $n$ and tends to
the "true" asymptotic value given by (\ref{segh}), i.e. $1.85575...$.
For instance, we get $c(5) \approx 1.73, \; \; c(6) \approx 1.79, \; \; c(7) \approx 1.82,$
and  $c(8) \approx 1.836.$
%************************* References ***********************************

\end{document}